\documentclass[default]{svmult}

\usepackage{mathptmx}       
\usepackage{helvet}         
\usepackage{courier}        
\usepackage{type1cm}        
\usepackage{makeidx}          
\usepackage{graphicx}         
\usepackage{multicol}         
\usepackage[bottom]{footmisc} 
\usepackage[english]{babel}                                        %
\usepackage[babel=true]{microtype}                                 %
\usepackage{amsmath}

\makeindex             

\begin{document}
\def\T{{\mbox{\rm\tiny T}}} 

\title*{Inexact Newton method for minimization of convex piecewise quadratic functions}
\titlerunning{Inexact Newton method}

\author{\bf A.I.~Golikov and I.E.~Kaporin}

\institute{%
A.I.~Golikov
\at Dorodnicyn Computing Center of FRC CSC RAS,
Vavilova 40, 119333 Moscow, Russia; 
\at Moscow Institute of Physics and Technology (State University), 
9 Institutskiy per., Dolgoprudny, Moscow Region, 141701, Russia; 
\email{gol-a@yandex.ru}
\and
I.E.~Kaporin
\at Dorodnicyn Computing Center of FRC CSC RAS,
Vavilova 40, 119333 Moscow, Russia; 
\email{igorkaporin@mail.ru}
}
\maketitle

\abstract{
An inexact Newton type method for numerical minimization of convex piecewise 
quadratic functions is considered and its convergence is analyzed. Earlier, 
a similar method was successfully applied to optimizaton problems arising in 
numerical grid generation. The method can be applied for computing a minimum 
norm nonnegative solution of underdetermined system of linear equations or 
for finding the distance between two convex polyhedra. The performance of the 
method is tested using sample data from NETLIB family of the University of 
Florida sparse matrix collection as well as quasirandom data.
} 
\section{Introduction} 
\label{intro}
The present paper is devoted to theoretical and experimental study of novel 
techniques for incorporation of preconditioned conjugate gradient linear 
solver into inexact Newton method. Earlier, similar method was successfully 
applied to optimizaton problems arising in numerical grid generation 
\cite{GK99,GKK04,Ka03}, and here we will consider its application to the 
numerical solution of piecewise-quadratic unconstrained optimization 
problems  \cite{GGEN09,KMPN17,Ma02,Ma04}. The latter include such problems 
as finding the projection of a given point onto the set of nonnegative 
solutions of an underdetermined system of linear equations \cite{GGE18} 
or finding a distance between two convex polyhedra \cite{Bo89} (and 
both are tightly related to the standard linear programming problem). 
The paper is organized as follows. In Section~2, a typical problem of 
minimization of piecewise quadratic function is formulated. In Section~3, 
certain technical results are given related to objective functions under 
consideration. Section~4 describes an inexact Newton method adjusted to 
the optimization problem. In Section~5, a convergence analysis of the proposed
algorithm is given with account of special stopping rule of inner linear 
conjugate gradient iterations. In Section~6, numerical results are presented
for various model problems. 

\section{Optimization problem setting}
\label{underdet}
Consider the piecewise-quadratic unconstrained optimization problem
\begin{equation}
p_* = \arg\min_{p\in R^m}
\left(\frac12\|(\widehat x + A^{\rm T}p)_+\|^2 - b^{\rm T}p\right),
\label{prob_set}
\end{equation}
where the standard notation $\xi_+=\max(0,\xi)=(\xi + |\xi|)/2$ is used.
Problem (\ref{prob_set}) can be viewed as the dual for finding  
projection of a vector on the set of nonnegative solutions of 
underdetermined linear systems of equations \cite{GGE18,GGEN09}:
$$
x_* = \arg\min_{^{Ax=b}_{x\ge 0}} \frac12 \|x-\widehat x\|^2,
$$
the solution of which is expressed via $p_*$ as 
$x_* = (\widehat x + A^{\rm T}p_*)_+$.
Therefore, we are considering piecewise quadratic 
function $\varphi:R^m\rightarrow R^1$ determined as
\begin{equation}
\varphi (p) = 
\frac12\|(\widehat x + A^{\rm T}p)_+\|^2 - b^{\rm T}p,
\label{phi}
\end{equation}
which is convex and differentiable.
Its gradient $g(p) = {\rm grad}~p$ is given by
\begin{equation}
g(p) = A(\widehat x + A^{\rm T}p)_+ - b,
\label{grad}
\end{equation}
and it has generalized Hessian \cite{HSN84}
\begin{equation}
H(p) = 
A{\rm Diag}\left({\rm sign}(\widehat x + A^{\rm T}p)_+\right)A^{\rm T}.
\label{genh}
\end{equation}
The relation of $H(p)$ to $\varphi(p)$ and $g(p)$ will be explained 
later in Remark~1.

%
\section{Taylor expansion of $(\cdot)_+^2$ function}
The following result is a special case of Taylor expansion 
with the residual term in integral form.

{LEMMA~1.} 
{\it
For any real scalars $\eta$ and $\zeta$ it holds
\begin{equation}
\frac12((\eta+\zeta)_+)^2 - \frac12(\eta_+)^2 - \zeta\eta_+
= \zeta^2\int_0^1\left(\int_0^1 {\rm sign}(\eta+st\zeta)_+ds\right)tdt.
\label{intexp1}
\end{equation}
}

{PROOF.}
Consider $f(\xi)=\frac12(\xi_+)^2$ and note that $f'(\xi)=\xi_+$ 
and $f''(\xi)={\rm sign}(\xi_+)$ (note that $f''(0)$ can formally be 
set equal to any finite real number, and w.l.o.g. we use $f''(0)=0$). 
Inserting this into the Taylor expansion
$$
f(\eta+\zeta) = f(\eta) + \zeta f'(\eta) +
\zeta^2 \int_0^1 \left(\int_0^1 f''(\eta+st\zeta)ds \right) tdt
$$
readily gives the desired result.

{LEMMA~2.} 
{\it
For any real $n$-vectors $y$ and $z$ it holds
\begin{equation}
\frac12\|(y+z)_+\|^2-\frac12\|y_+\|^2-z^{\rm T}y_+ 
= \frac12 z^{\rm T}{\rm Diag}(d)z,
\label{intexp1a}
\end{equation}
where
\begin{equation}
d = \int_0^1\left(\int_0^1 {\rm sign}(y+stz)_+ds\right)2tdt.
\label{intexp1b}
\end{equation}
}

{PROOF.}
Setting in (\ref{intexp1}) $\eta=y_j$, $\zeta = z_j$, 
and summing over all $j=1,\ldots,n$ obviously yields 
the required formula. Note that the use of scalar 
multiple 2 within the integral provides for the estimate
$\|{\rm Diag}(d)\| \le 1$.

{LEMMA~3.} 
{\it 
Function (\ref{phi}) and its gradient (\ref{grad}) satisfy the identity 
\begin{equation}
\varphi(p+q)-\varphi(p)-q^{\rm T}g(p) 
= \frac12 q^{\rm T}A\,{\rm Diag}(d)A^{\rm T}q,
\label{intexp2a}
\end{equation}
where
\begin{equation}
d = \int_0^1\left(\int_0^1 
{\rm sign}(\widehat x + A^{\rm T}p + stA^{\rm T}q)_+ds\right)2tdt.
\label{intexp2b}
\end{equation}
}

{PROOF.}
Setting in (\ref{intexp1a}) and (\ref{intexp1b}) $y=\widehat x + A^{\rm T}p$ 
and $z = A^{\rm T}q$ readily yields the required result (with account 
of cancellation of linear terms involving $b$ in the left hand side of 
(\ref{intexp2a})).

{REMARK~1.} As is seen from (\ref{intexp2b}), if the condition 
\begin{equation}
{\rm sign}(\widehat x + A^{\rm T}p + \vartheta A^{\rm T}q)_+ 
= {\rm sign}(\widehat x + A^{\rm T}p)_+, 
\label{nearsol}
\end{equation}
holds true for any $0\le\vartheta\le 1$, then (\ref{intexp2a}) is simplified as
\begin{equation}
\varphi(p+q) - \varphi(p) - q^{\rm T}g(p) = \frac12 q^{\rm T}H(p)q,
\label{locquad}
\end{equation}
where the generalized Hessian matrix $H(p)$ is defined 
in (\ref{genh}). This explains the key role of $H(p)$ in 
the organization of the Newton-type method considered below. 
Note that a sufficient condition for (\ref{nearsol}) to hold is 
\begin{equation}
|(A^{\rm T}q)_j| \le |(\widehat x + A^{\rm T}p)_j| \qquad
{\rm whenever} \qquad (A^{\rm T}q)_j (\widehat x + A^{\rm T}p)_j < 0;
\label{nearsol2}
\end{equation}
that is, if certain components of the increment $q$ are relatively small, 
then $\varphi$ is exactly quadratic (\ref{locquad}) in the corresponding 
neighborhood of $p$.

%
\section{Inexact Newton method for dual problem}

As suggests condition (\ref{nearsol}) and its consequence 
(\ref{locquad}), one can try to numerically minimize $\varphi$ 
using Newton type method $p_{k+1} = p_k - d_k$, where 
$d_k = H(p_k)^{-1}g(p_k)$.
Note that by (\ref{locquad}) this will immediately give 
the exact minimizer $p_* = p_{k+1}$ if the magnitudes 
of $d_k$ components are sufficiently small to satisfy 
(\ref{nearsol2}) taken with $p=p_k$ and $q=-d_k$.
However, initially $p_k$ may be rather far from solution, 
and only gradual improvements are possible. First, a damping 
factor $\alpha_k$ must be used to guarantee monotone 
convergence (with respect to the decrease of $\varphi(p_k)$ 
as $k$ increases).  Second, $H(p_k)$ must be replaced 
by some appropriate approximation $M_k$ in order to provide its
invertibility with a reasonable bound for the inverse. Therefore, 
we propose the following prototype scheme 
$$
p_{k+1} = p_k - \alpha_k M_k^{-1}g(p_k), 
$$
where
\begin{equation}
M_k = H(p_k) + \delta {\rm Diag}(AA^{\rm T}).
\label{reg_hess}
\end{equation}
The parameters $0< \alpha_k \le 1$ and $0\le \delta \ll 1$ must 
be defined properly for better convergence. Furthermore, at 
initial stages of iteration, the most efficient strategy is 
to use approximate Newton directions $d_k \approx M_k^{-1}g(p_k)$,
which can be obtained using preconditioned  conjugate gradient (PCG) 
method for the solution of Newton equation $M_kd_k=g(p_k)$. As will be 
seen later, it suffices to use any vector $d_k$ which satisfies 
conditions
\begin{equation}
d_k^{\rm T}g_k = d_k^{\rm T}M_kd_k =\vartheta_k^2g_k^{\rm T}M_k^{-1}g_k
\label{newt_dir}
\end{equation}
with \mbox{$0<\vartheta_k<1$} sufficiently separated from zero. 
For any preconditioning, the approximations constructed by the PCG 
method satisfy (\ref{newt_dir}), see Section~{\ref{PCG}} below. 
With account of the Armijo type criterion 
\begin{equation}
\varphi(p_k - \alpha d_k) \le 
\varphi(p_k) - \frac{\alpha}{2} d_k^{\rm T}g(p_k),
\qquad \alpha\in\lbrace{1,~1/2,~1/4,~\ldots\rbrace},
\label{back_trk}
\end{equation}
where the maximum steplength $\alpha$ satisfying (\ref{back_trk}) is used, 
the inexact Newton algorithm can be presented as follows:

\medskip
{\bf Algorithm~1.}

{\it Input:}

$A \in R^{m\times n}$, \quad 
$b\in R^m$,  \quad 
$\widehat x\in R^m$;

{\it Initialization:}

$\delta = 10^{-6}$, \quad 
$\varepsilon = 10^{-12}$, \quad
$\tau = 10^{-15}$

$k_{\max} = 2000$,  \quad $l_{\max} = 10$; \quad $p_0 = 0$, 

{\it Iterations:}

{\bf for} $k = 0, 1, \ldots, k_{\max}-1$:

~~~~~~$x_k = (\widehat x + A^{\rm T}p_k)_+$

~~~~~~$\varphi_k = \frac12 \|x_k\|^2 - b^{\rm T}p_k$

~~~~~~$g_k = Ax_k - b$

~~~~~~{\bf if} $(\|g_k\| \le \varepsilon\|b\|)$ 
{\bf return} $\lbrace{x_k,p_k,g_k\rbrace}$ 
%

~~~~~~{\bf find} $d_k \in R^m$ ~{\bf such that}

~~~~~~~~~~~~$d_k^{\rm T}g_k = d_k^{\rm T}M_kd_k = \vartheta_k^2g_kM_k^{-1}g_k$,

~~~~~~~~~~~~{\bf where} $M_k = A~{\rm Diag}({\rm sign}(x_k))~A^{\rm T} 
+ \delta {\rm Diag}(AA^{\rm T})$

~~~~~~$\alpha^{(0)} = 1$

~~~~~~$p_k^{(0)} = p_k - d_k$

~~~~~~{\bf for} $l = 0, 1, \ldots, l_{\max}-1$:

~~~~~~~~~~~~$x_k^{(l)} = (\widehat x + A^{\rm T}p_k^{(l)})_+$

~~~~~~~~~~~~$\varphi_k^{(l)} = \frac12 \|x_k^{(l)}\|^2 - b^{\rm T}p_k^{(l)}$

~~~~~~~~~~~~$\zeta_k^{(l)} = \left(\frac12
\alpha^{(l)}d_k^{\rm T}g_k + \varphi_k^{(l)}\right)-\varphi_k$

~~~~~~~~~~~~{\bf if} ($\zeta_k^{(l)}>\tau|\varphi_k|$) {\bf then}

~~~~~~~~~~~~~~~~~~$\alpha^{(l+1)} = \alpha^{(l)}/2$

~~~~~~~~~~~~~~~~~~$p_k^{(l+1)} = p_k - \alpha^{(l+1)} d_k$

~~~~~~~~~~~~{\bf else}

~~~~~~~~~~~~~~~~~~$p_k^{(l+1)} = p_k^{(l)}$ 

~~~~~~~~~~~~~~~~~~{\bf go to} NEXT

~~~~~~~~~~~~{\bf end if}

~~~~~~{\bf end for}

~~~~~~NEXT:~$p_{k+1} = p_k^{(l+1)}$

{\bf end for}
%

\medskip
Next we explore the convergence properties of this algorithm.

%
\section{Convergence analysis of inexact Newton method}
It appears that Algorithm~1 exactly conforms with the convergence
analysis presented in \cite{Ka03} (see also \cite{AK01}). 
For the completeness of presentation and compatibility of notations, 
we reproduce here the main results of \cite{Ka03}. 

\subsection{Estimating convergence of inexact Newton method}
The main assumptions we need for the function $\varphi(p)$ 
under consideration are that it is bounded from below, have 
gradient $g(p)\in R^m$, and satisfies
\begin{equation}
\varphi(p+q)-\varphi(p)-q^{\rm T}g(p) \le {\gamma\over 2}q^{\rm T}Mq
\label{uppbnd}
\end{equation}
for the symmetric positive definite $m\times m$ matrix $M=M(p)$ defined 
above in (\ref{reg_hess}) and some constant $\gamma\ge 1$.
Note that the exact knowledge of $\gamma$ is not necessary for actual 
calculations. The existence of such $\gamma$ follows from (\ref{reg_hess}) 
and Lemma~3. Indeed, denoting 
$D = ({\rm Diag}(AA^{\rm T}))^{1/2}$ and $\widehat A = D^{-1}A$, 
for the right hand side of (\ref{intexp2a}) one has, with account 
of $\|{\rm Diag}(d)\|\le 1$ and $H(p)\ge 0$,
$$
A{\rm Diag}(d)A^{\rm T} \le AA^{\rm T} 
\le \|\widehat A\|^2{\rm Diag}(AA^{\rm T}) 
\le \frac{\|\widehat A\|^2}{\delta} 
\Bigl(\delta{\rm Diag}(AA^{\rm T})+H(p)\Bigr)
=  \frac{\|\widehat A\|^2}{\delta} M; 
$$
therefore, (\ref{uppbnd}) holds with 
\begin{equation}
\gamma = \|\widehat A\|^2/\delta.
\label{gamma}
\end{equation}
The latter formula explains our choice of $M$ which is more appropriate
in cases of large variations in norms of rows in $A$ (see the examples 
and discussion in Section~\ref {NumTest1}). 


Next we will estimate the reduction in the value of $\varphi$ attained
by the descent along the direction $(-d)$ satisfying  (\ref{newt_dir}). 
One can show the following estimate for the decrease of objective function 
value at each iteration (here, simplified notations $p=p_k$, $\hat p = p_{k+1}$ 
etc. are used) \mbox{$\hat p = p - \alpha d$} with $\alpha=2^{-l}$, 
where $l=0,1,\ldots$, as evaluated according 
to (\ref{back_trk}):
\begin{equation}
\varphi(\hat p) \le \varphi(p) - \frac{\vartheta^2}{4\gamma} g^{\rm T}M^{-1}g.
\label{g_conv}
\end{equation}
In particular, if the values of $\vartheta^2$ are separated from zero by 
a positive constant $\vartheta_{\min}^2$ (lower estimate for $\vartheta$
follows from Section~\ref{PCG} and an upper bound for 
\mbox{$\kappa={\rm cond}(CM)$)}, then, with account for 
$M \le (1+\delta)\|A\|^2 I$ and the boundedness of $\varphi$ from below, 
it follows 
$$
\sum_{j=0}^{k-1}g_j^{\rm T}g_j 
\le \frac{4\gamma(1+\delta)\|A\|^2
(\varphi(p_0)-\varphi(p_*))}{\vartheta_{\min}^2}.
$$
Noting that the right hand side of the latter estimate does not depend 
on $k$, it finally follows that
$$
\lim_{k\rightarrow\infty}\|g(p_k)\| = 0,
$$
where $k$ is the number of the outer (nonlinear) iteration.

Estimate (\ref{g_conv}) can be verified as follows (note that quite 
similar analysis can be found in \cite{Ma02}). Using $q = -\beta d$, 
where $0 < \beta < 2/\gamma$, one can obtain from (\ref{uppbnd}) and 
(\ref{newt_dir}) the following estimate for the decrease of $\varphi$ 
along the direction $(-d)$:
\begin{eqnarray}
\nonumber
\varphi(p-\beta d) 
&=& \varphi(p) - \beta d^{\rm T}g +
\left(\varphi(p-\beta d) - \varphi(p) - \beta d^{\rm T}g\right) \\
\label{upp2bnd}
&\le& \varphi(p) -
\left(\beta - {\gamma\over 2}\beta^2\right)d^{\rm T}Md.\\
\nonumber
\end{eqnarray}
The following two cases are possible.

{\it Case~1.} If the condition (\ref{back_trk}) is satisfied at once 
for $\alpha=1$, this means that (recall that the left equality of 
(\ref{newt_dir}) holds)
$$
\varphi(\hat p) \le \varphi(p) - \frac12 d^{\rm T}Md.
$$

{\it Case~2.} Otherwise, if at least one bisection of steplength 
was performed (and the actual steplength is $\alpha$), then, using 
(\ref{upp2bnd}) with $\beta = 2\alpha$, it follows
$$
\varphi(p) - \alpha d^{\rm T}Md < \varphi(p - 2\alpha d)
\le \varphi(p) - (2\alpha - 2\gamma\alpha^2) d^{\rm T}Md,
$$
which readily yields $\alpha > 1/(2\gamma)$. Since we also have
$$
\varphi(p - \alpha d) \le \varphi(p) - {\alpha\over 2} d^{\rm T}Md,
$$
it follows
\begin{equation}
\varphi(p - \alpha d) \le \varphi(p) - {1\over 4\gamma} d^{\rm T}Md.
\label{g2conv}
\end{equation}
Joining these two cases, taking into account that $\gamma\ge 1$, and
using the second equality in (\ref{newt_dir}) one obtains the required 
estimate (\ref{g_conv}).

It remains to notice that as soon as the norms of $g$ attain sufficiently 
small values, the resulting directions $d$ will also have small norms. 
Therefore, the case considered in Remark~1 will take place, and 
finally the convergence of the Newton method will be much faster than 
at its initial stage.

\subsection{Linear CG approximation of Newton directions}
\label{stopCG}
Next we relate the convergence of inner linear Preconditioned 
Conjugate Gradient (PCG) iterations to the efficiency of 
Inexact Newton nonlinear solver. Similar issues were considered
in \cite{AK01,GKK04,KA94,Ka03}.

An approximation $d^{(i)}$ to the solution of the problem $Md=g$ generated 
on the $i$th PCG iteration by the recurrence $d^{(i+1)} = d^{(i)} + s^{(i)}$ 
(see Algorithm~2 below) can be written as follows (for our purposes, 
we always set the initial guess for the solution $d^{(0)}$ to zero): 
\begin{equation}
d^{(i)} = \sum_{j=0}^{i-1}s^{(j)},
\label{pcg_di}
\end{equation}
where the PCG direction vectors are pairwise $M$-orthogonal:
\mbox{$(s^{(j)})^{\T}Ms^{(l)}=0$}, ~\mbox{$j\ne l$}.
Let also denote the $M$-norms of PCG directions as
$\eta^{(j)} = (s^{(j)})^{\T}Ms^{(j)}$, $j=0,1,\ldots,i-1$.
Therefore, from (\ref{pcg_di}), one can determine  
$$
\zeta^{(i)} = (d^{(i)})^{\T} M d^{(i)} = \sum_{j=0}^{i-1}\eta^{(j)},
$$
and estimate (\ref{g2conv}) takes the form 
$$
\varphi_{k+1} \le \varphi_k - {1\over 4\gamma} \sum_{j=0}^{i_k-1}\eta_k^{(j)}, 
$$
where $k$ is the Newton iteration number. Summing up the latter
inequalities for $0 \le k \le m-1$, we get
\begin{equation}
c_0 \equiv 4\gamma(\varphi_0 -  \varphi_*)  
\ge  \sum_{k=0}^{m-1} \sum_{j=0}^{i_k-1}\eta_k^{(j)}
\label{est_conv}
\end{equation}
On the other hand, the cost measure related to the total time needed
to perform $m$ inexact Newton iterations with $i_k$ PCG iterations at 
each Newton step, can be estimated as proportional to
$$
T_m=\sum_{k=0}^{m-1} \left( \epsilon_{\rm CG}^{-1} + i_k\right)
\le c_0\frac{\sum_{k=0}^{m-1}\left( \epsilon_{\rm CG}^{-1} + i_k\right)}
{\sum_{k=0}^{m-1} \sum_{j=0}^{i_k-1}\eta_k^{(j)}}
\le c_0 \max_{k<m}\frac{\epsilon_{\rm CG}^{-1} + i_k}{\sum_{j=0}^{i_k-1}\eta_k^{(j)}}. 
$$
Here $\epsilon_{\rm CG}$ is a small parameter reflecting the ratio of 
one linear PCG iteration cost to the cost of one Newton iteration 
(in particular, including construction of preconditioning and 
several $\varphi$ evaluations needed for backtracking) plus 
possible efficiency loss due to early PCG termination.
Thus, introducing the function 
\mbox{$\psi(i) = (\epsilon_{\rm CG}^{-1} + i)/\zeta^{(i)}$},
(here, we omit the index $k$) one obtains a reasonable criterion 
to stop PCG iterations in the form \mbox{$\psi(i) > \psi(i-1)$}.
Here, the use of smaller values $\varepsilon_{\rm CG}$ generally 
corresponds to the increase of the resulting iteration number bound.
Rewriting the latter condition, one obtains the final form of the 
PCG stopping rule: 
\begin{equation}
(\epsilon_{\rm CG}^{-1} + i)\eta^{(i-1)} \le \zeta^{(i)}.
\label{new_stop}
\end{equation}
Note that by this rule, the PCG iteration number is always no less than 2.

Finally, we explicitly present the resulting  formulae for the 
PCG algorithm incorporating the new stopping rule.
Following \cite{GGEN09}, we use the Jacobi preconditioning 
\begin{equation}
C = ({\rm Diag}(M))^{-1}.
\label{diag_prec}
\end{equation}
Moreover, the reformulation \cite{KM11} of the CG algorithm \cite{HS52,Ax76} 
is used. This may give a more efficient parallel implementation, see, e.g., 
\cite{GGEN09}.

Following \cite{KM11}, recall that at each PCG iteration the 
$M^{-1}$-norm of the \mbox{$(i+1)$}-th residual \mbox{$r^{(i+1)}=g-Md^{(i+1)}$} 
attains its minimum over the corresponding Krylov subspace. 
Using the standard PCG recurrences (see Section~\ref{PCG} below) 
one can find 
$d^{(i+1)} = d^{(i)}+Cr^{(i)}\alpha^{(i)}+s^{(i-1)}\alpha^{(i)}\beta^{(i-1)}$.
Therefore, the optimum increment $s^{(i)}$ in the recurrence 
\mbox{$d^{(i+1)} = d^{(i)} + s^{(i)}$}, where \mbox{$s^{(i)} = V^{(i)} h^{(i)}$} 
and \mbox{$V^{(i)} = [ Cr^{(i)} ~|~ s^{(i-1)}]$}, can be determined via 
the solution of the following 2-dimensional linear least squares problem:
$$
\left[_{\beta^{(i)}}^{\alpha^{(i)}}\right] = h^{(i)}
=\arg\min_{h\in{\bf R}^2}\|g-Md^{(i+1)}\|_{M^{-1}}
=\arg\min_{h\in{\bf R}^2}\|r^{(i)}-MV^{(i)} h\|_{M^{-1}}. 
$$
By redefining \mbox{$r^{(i)}:=-r^{(i)}$} and introducing vectors 
\mbox{$t^{(i)}=Ms^{(i)}$}, the required PCG reformulation follows:
\bigskip

{\bf Algorithm~2.}

$r^{(0)} = -g,$\quad 
$d^{(0)} = s^{(-1)} = t^{(-1)} = 0,$ \quad 
$\zeta^{(-1)} = 0;$

$i=0,1,\ldots,it_{\max}:$

\qquad $w^{(i)} = Cr^{(i)},$

\qquad $z^{(i)} = Mw^{(i)},$

\qquad $\gamma^{(i)} = (r^{(i)})^{\T}w^{(i)},$ 
\quad $\xi^{(i)} = (w^{(i)})^{\T}z^{(i)},$ 
\quad $\eta^{(i-1)} = (s^{(i-1)})^{\T}t^{(i-1)},$

\qquad $\zeta^{(i)} = \zeta^{(i-1)} + \eta^{(i-1)},$

\qquad {\bf if}
~$((\varepsilon_{\rm CG}^{-1}+i)\eta^{(i-1)} \le\zeta^{(i)})$~ 
{\bf or} ~$(\gamma^{(i)} \le \varepsilon_{\rm CG}^2  \gamma^{(0)})$~
{\bf return}~$\lbrace d^{(i)}\rbrace$;

\qquad {\bf if} ~ ($k=0$) ~ {\bf then} 

\qquad\qquad $\alpha^{(i)} = -\gamma^{(i)}/\xi^{(i)}, \quad \beta^{(i)} = 0;$

\qquad {\bf else} 

\qquad\qquad \mbox{$\delta^{(i)}
=\gamma^{(i)}/(\xi^{(i)}\eta^{(i-1)}-(\gamma^{(i)})^2)$}, 

\qquad\qquad \mbox{$\alpha^{(i)} = -\eta^{(i-1)}\delta^{(i)}$}, 
\quad\qquad \mbox{$\beta^{(i)} = \gamma^{(i)}\delta^{(i)}$};

\qquad {\bf end if} 

\qquad $t^{(i)} = z^{(i)}\alpha^{(i)} + t^{(i-1)}\beta^{(i)},$~ 
\qquad $r^{(i+1)} = r^{(i)} + t^{(i)},$

\qquad $s^{(i)} = w^{(i)}\alpha^{(i)} + s^{(i-1)}\beta^{(i)},$ 
\qquad $d^{(i+1)} = d^{(i)} + s^{(i)}.$

\bigskip
\noindent
For maximum reliability, the new stopping rule (\ref{new_stop}) 
is used along with the standard one; however, in almost all cases 
the new rule provides for an earlier CG termination.

Despite of somewhat larger workspace and number of vector operations 
compared to the standard algorithm, the above version of CG algorithm 
enables more efficient parallel implementation of scalar product operations. 
At each iteration of the above presented algorithm, it suffices to 
use one {\small MPI$\_$AllReduce($*$,$*$,3,\dots)}
operation instead of two {\small MPI$\_$AllReduce($*$,$*$,1,\dots)}
operation in the standard PCG recurrences. This is especially important 
when many MPI processes are used and the start-up time for MPI$\_$AllReduce operations is relatively large.
For another equivalent PCG reformulations allowing to properly reorder the 
scalar product operations, see \cite{DE03} and references cites therein.

%
\subsection{Convergence properties of PCG iterations}
\label{PCG}
Let us recall some basic properties of the PCG algorithm, see, e.g. 
\cite{Ax76}. The standard PCG algorithm (algebraically equivalent to
Algorithm~2) for the solution of the problem $Md=g$ can be written 
as follows (the initial guess for the solution $d_0$ is set to zero): 
\begin{eqnarray}
&&d^{(0)} = 0, \quad r^{(0)} = g, \quad s^{(0)} = Cr^{(0)}; 
\nonumber
\\
&&{\bf for} ~ i = 0, 1, \ldots, m-1:
\nonumber
\\
&&\qquad \alpha^{(i)} = (r^{(i)})^TCr^{(i)}/(s^{(i)})^TMs^{(i)},
\nonumber
\\
&&\qquad d^{(i+1)} = d^{(i)} + s^{(i)}\alpha^{(i)},
\nonumber
\\
&&\qquad r^{(i+1)} = r^{(i)} - Ms^{(i)}\alpha^{(0)},
\nonumber
\\
&&\qquad {\bf if} ~ ((r^{(i)})^TCr^{(i+1)} \le 
\varepsilon_{\rm CG}^2 (r^{(0)})^TCr^{(0)}) ~ {\bf return} ~ d^{(i+1)}
\nonumber
\\
&&\qquad \beta^{(i)} = (r^{(i+1)})^TCr^{(i+1)}/(r^{(i)})^TCr^{(i)},
\nonumber
\\
&&\qquad s^{(i+1)} = Cr^{(i+1)} + s^{(i)}\beta^{(i)}.
\\
&&{\bf end for}
\nonumber
\end{eqnarray}
The scaling property (\ref{newt_dir}) (omitting the upper and lower 
indices at $d$, it reads $d^{\rm T}g = d^{\rm T}Md$) can be proved 
as follows. Let $d=d^{(i)}$ be obtained after $i$ iterations of the PCG 
method applied to $Md=g$ with zero initial guess $d^{(0)}=0$. Therefore,
$d \in K_i = {\rm span} \lbrace Cg, CMCg, \dots, (CM)^{i-1}Cg\rbrace$,
and, by the PCG optimality property, it holds
$$
d = \arg\min_{d\in K_i} (g-Md)^TM^{-1}(g-Md).
$$
Since $\alpha d\in K_i$ for any scalar $\alpha$, one gets
$$
(g-\alpha Md)^TM^{-1}(g-\alpha Md) \ge (g-Md)^TM^{-1}(g-Md).
$$
Setting here $\alpha=d^Tg/d^TMd$, one can easily transform this 
inequality as $0 \ge (-d^Tg+d^TMd)^2$, which readily yields 
(\ref{newt_dir}). Furthermore, by the well known estimate of the 
PCG iteration error \cite{Ax76} using Chebyshev polynomials, one gets
$$
1-\theta^2 \equiv (g-Md)^TM^{-1}(g-Md)/g^TM^{-1}g \le
\cosh^{-2}\left(2i/\sqrt{\kappa}\right)
$$
where
$$
\kappa = {\rm cond}(CM) \equiv \lambda_{\max}(CM)/\lambda_{\min}(CM).
$$
By the scaling condition, this gives
\begin{equation}
\theta^2 = d^TMd/g^TM^{-1}g \ge \tanh^2\left(2i/\sqrt{\kappa}\right).
\label{pcg_conv}
\end{equation}
Hence, $0<\theta<1$ and $\theta^2\rightarrow 1$ 
as the PCG iteration number $i$ grows.

\section{Numerical test results}
\label{NumTest1}
Below we consider two families of test problems which can be solved
via minimization of piecewise quadratic problems. The first one was 
described above in Section~\ref{underdet} (see also \cite{GGE18}), 
while the second coincides with the problem setting for the evaluation 
of distance between two convex polyhedra used in \cite{Bo89}.
The latter problem is of key importance e.g., in robotics and computer
animation.
 
\subsection{Test results for 11 NETLIB problems}
Matrix data from the following 11 linear programming problems 
(this is the same selection from NETLIB collection as considered 
in \cite{KMPN17}), were used to form test problems (\ref{prob_set}).
Note that further we only consider the case $\widehat x=0$.
Recall also the notation $x_* = (\widehat x + A^Tp_*)_+$.
The problems in Table~1 below are ordered by the number of 
nonzero elements ${\rm nz}(A)$ in $A\in R^{m\times n}$.

\begin{table}
\begin{small}
\caption{Matrix properties for 11 NetLib problems}
\label{tab1}
\begin{tabular}{p{2cm}p{1cm}p{1cm}p{1cm}p{2cm}p{2cm}p{2cm}}
\hline\noalign{\smallskip}
name & m & n & nz(A) & $\|x_*\|$ & $\min_i(AA^T)_{ii}$ & $\max_i(AA^T)_{ii}$\\
\noalign{\smallskip}\hline\noalign{\smallskip}
afiro      &    27 &    51 &    102 & 634.029569 & 1.18490000 & 44.9562810\\
addlittle  &    56 &   138 &    424 & 430.764399 & 1.00000000 & 10654.0000\\
agg3       &   516 &   758 &   4756 & 765883.022 & 1.00000001 & 179783.783\\
25fv47     &   821 &  1876 &  10705 & 3310.45652 & 0.00000000 & 88184.0358\\
pds$\_$02  &  2953 &  7716 &  16571 & 160697.180 & 1.00000000 & 91.0000000\\
cre$\_$a   &  3516 &  7248 &  18168 & 1162.32987 & 0.00000000 & 27476.8400\\
80bau3b    &  2262 & 12061 &  23264 & 4129.96530 & 1.00000000 & 321739.679\\
ken$\_$13  & 28362 & 42659 &  97246 & 25363.3224 & 1.00000000 & 170.000000\\
maros$\_$r7&  3136 &  9408 & 144848 & 141313.207 & 3.05175947 & 3.37132546\\
cre$\_$b   &  9648 & 77137 & 260785 & 624.270129 & 0.00000000 & 27476.8400\\
osa$\_$14  &  2337 & 54797 & 317097 & 119582.321 & 18.0000000 & 845289.908\\
\noalign{\smallskip}\hline\noalign{\smallskip}
\end{tabular}
\end{small}
\end{table}

\begin{table}
\begin{small}
\caption{Computational results of different solvers for 11 NetLib problems}
\label{tab2}
\begin{tabular}{p{2cm}p{2cm}p{2cm}p{1.7cm}p{1.7cm}p{1.7cm}}
\hline\noalign{\smallskip}
name & solver & time(sec) & $\|Ax - b\|_\infty$ & {\#}NewtIter & {\#}MVMult \\
\noalign{\smallskip}\hline\noalign{\smallskip}
afiro      & GNewtEGK  &  0.001 & 8.63E--11 &  17 &    398 \\
--"--      & ssGNewton &  0.06  & 6.39E--14 &  -- &     -- \\
--"--      & cqpMOSEK  &  0.31  & 1.13E--13 &  -- &     -- \\
\hline\noalign{\smallskip}
addlittle  & GNewtEGK  &  0.003 & 6.45E--10 &  22 &   1050 \\
--"--      & ssGNewton &  0.05  & 2.27E--13 &  -- &     -- \\
--"--      & cqpMOSEK  &  0.35  & 7.18E--11 &  -- &     -- \\
\hline\noalign{\smallskip}
agg3       & GNewtEGK  &  0.14  & 3.93E--07 & 116 &   9234 \\
--"--      & ssGNewton &  0.27  & 3.59E--08 &  -- &     -- \\
--"--      & cqpMOSEK  &  0.40  & 2.32E--10 &  -- &     -- \\
\hline\noalign{\smallskip}
25fv47     & GNewtEGK  &  0.54  & 7.15E--10 & 114 &  32234 \\
--"--      & ssGNewton &  1.51  & 3.43E--09 &  -- &     -- \\
--"--      & cqpMOSEK  &  1.36  & 1.91E--11 &  -- &     -- \\
\hline\noalign{\smallskip}
pds$\_$02  & GNewtEGK  &  0.32  & 1.55E--08 &  75 &   8559 \\
--"--      & ssGNewton &  2.30  & 1.40E--07 &  -- &     -- \\
--"--      & cqpMOSEK  &  0.51  & 8.20E--06 &  -- &     -- \\
\hline\noalign{\smallskip}
cre$\_$a   & GNewtEGK  &  3.36  & 2.64E--09 & 219 &  85737 \\
--"--      & ssGNewton &  1.25  & 4.13E--06 &  -- &     -- \\
--"--      & cqpMOSEK  &  0.61  & 2.15E--10 &  -- &     -- \\
\hline\noalign{\smallskip}
80bau3b    & GNewtEGK  &  0.27  & 3.33E--09 &  79 &   6035 \\
--"--      & ssGNewton &  0.95  & 1.18E--12 &  -- &     -- \\
--"--      & cqpMOSEK  &  0.80  & 2.90E--07 &  -- &     -- \\
\hline\noalign{\smallskip}
ken$\_$13  & GNewtEGK  &  1.41  & 2.70E--08 &  55 &   6285 \\
--"--      & ssGNewton &  9.09  & 4.39E--09 &  -- &     -- \\
--"--      & cqpMOSEK  &  2.09  & 1.71E--09 &  -- &     -- \\
\hline\noalign{\smallskip}
maros$\_$r7& GNewtEGK  &  0.10  & 1.18E--09 &  27 &    535 \\
--"--      & ssGNewton &  2.86  & 2.54E--11 &  -- &     -- \\
--"--      & cqpMOSEK  & 55.20  & 3.27E--11 &  -- &     -- \\
\hline\noalign{\smallskip}
cre$\_$b   & GNewtEGK  &  9.25  & 6.66E--10 &  75 &  24590 \\
--"--      & ssGNewton & 13.20  & 1.62E--09 &  -- &     -- \\
--"--      & cqpMOSEK  &  2.31  & 1.61E--06 &  -- &     -- \\
\hline\noalign{\smallskip}
osa$\_$14  & GNewtEGK  & 42.59  & 8.25E--08 & 767 & 104874 \\
--"--      & ssGNewton & 60.10  & 4.10E--08 &  -- &     -- \\
--"--      & cqpMOSEK  &  4.40  & 7.82E--05 &  -- &     -- \\
\noalign{\smallskip}\hline\noalign{\smallskip}
\end{tabular}
\end{small}
\end{table}

\begin{table}
\begin{small}
\caption{Comparison of JCG stopping criteria 
$\gamma^{(i)}\le\varepsilon_{CG}^2\gamma^{(0)}$ (old) 
and (\ref{new_stop}) (new)} 
\label{tab3}
\begin{tabular}{p{2.2cm}p{2.2cm}p{2.1cm}p{2.2cm}p{2.2cm}}
\hline\noalign{\smallskip}
criterion &$\varepsilon_{CG}$& geom.mean time & arithm.mean time 
& geom.mean res.\\
\noalign{\smallskip}\hline\noalign{\smallskip}
old      &    0.05 &   1.78 &  12.37 & 3.64--09 \\
old      &    0.03 &   1.45 &  12.21 & 3.47--09 \\
old      &    0.01 &   1.35 &  14.91 & 2.11--09 \\
old      &   0.003 &   1.48 &  21.14 & 2.14--09 \\
old      &   0.001 &   1.82 &  29.12 & 3.64--09 \\
\noalign{\smallskip}\hline\noalign{\smallskip}
new      &   0.003 &   1.66 &  11.12 & 3.46--09 \\
new      &   0.002 &   1.39 &  11.36 & 3.64--09 \\
new      &   0.001 &   1.26 &  10.81 & 4.65--09 \\
new      &  0.0003 &   1.26 &  11.25 & 4.23--09 \\
new      &  0.0001 &   1.33 &  12.47 & 3.60--09 \\
\noalign{\smallskip}\hline\noalign{\smallskip}
\end{tabular}
\end{small}
\end{table}

It is readily seen that 3 out of 11 matrices have null rows, and more than
half of them have rather large variance of row norms. This explains the 
proposed Hessian regularization (\ref{reg_hess}) instead of the earlier 
construction \cite{GGE18,KMPN17} $M_k = H(p_k) + \delta I_m$. The latter 
is a proper choice only for matrices with rows of nearly equal length, 
such as ~{\bf maros$\_$r7}~ example or various matrices with uniformly 
distributed quasirandom entries, as used for testing in 
\cite{GGEN09,KMPN17}. In particular, estimate (\ref{gamma}) with $D=I$ 
would take the form $\gamma=\|A\|^2/\delta$, so the resulting method appears 
to be rather sensitive to the choice of $\delta$.


In Table~2, the results presented in \cite{KMPN17} are reproduced 
along with similar data obtained with our version of Generalized 
Newton method. It must be stressed that we used the fixed set of 
tunung parameters 
\begin{equation}
\delta=10^{-6}, \qquad 
\varepsilon=10^{-12}, \qquad 
\varepsilon_{\rm CG}=10^{-3}, \qquad 
l_{\max} = 10,
\label{def_set}
\end{equation}
for all problems. Note that In \cite{KMPN17} the parameter choice 
for the Armijo procedure was not specified. 

In \cite{KMPN17}, the calculations were performed on 
5GHz AMD 64 Athlon X2 Dual Core. In our experiments, one core 
of 3.40 GHz x8 Intel (R) Core (TM) i7-3770 CPU was used, which 
is likely somewhat slower.

Note that the algorithm of \cite{KMPN17} is based on direct evaluation 
of $M_k$ and its sparse Cholesky factorization, while our 
implementation, as was proposed in \cite{GGEN09}, uses the Jacobi 
preconditioned Conjugate Gradient iterations for approximate 
evaluation of Newton directions. Thus, the efficiency of our 
implementation critically depends on the CG iteration convergence, 
which is sometimes slow. 
On the other hand, since the main computational kernels of the 
algorithm are presented by matrix-vector multiplications of the 
type ~$x=Ap$~ or ~$q=A^Ty$, its parallel implementation can be 
sufficiently efficient. 

In Table~2, the abbreviation {\bf cqpMOSEK} refers to MOSEK Optimization 
Software package for convex quadratic problems, see \cite{KMPN17}. 
The abbreviation {\bf ssGNewton} denotes the method implemented and 
tested in \cite{KMPN17}, while {\bf GNewtEGK} stands for the 
method proposed in the present paper. 

Despite the use of slower computer, our {\bf GNewtEGK} demonstrates 
considerably faster performance in 8 cases of 11. Otherwise, one can 
observe that smaller computational time of {\bf cqpMOSEK} goes along 
with much worse residual norm, see the results for problems 
~{\bf cre$\_$b}~ and ~{\bf osa$\_$14}~. 

Thus, in most cases the presented implementation of Generalized 
Newton method takes not too large number of Newton iterations 
using approximate Newton directions generated by CG iterations 
with diagonal preconditioning (\ref{diag_prec}) and special 
stopping rule (\ref{new_stop}). 

A direct comparison of efficiency for the standard CG iterations stopping 
rule \mbox{$\gamma^{(i)}\le\varepsilon_{\rm CG}^2\gamma^{(0)}$} 
(see Algorithm~2 for the notations) and the new one (\ref{new_stop}) 
is given in Table~\ref{tab3}, where the timing (in seconds) and 
precision results averaged over the same 11 problems are given. 
One can see that nearly the same average residual norm $\|Ax-b\|_{\infty}$ 
can be obtained considerably faster and with less critical dependence on 
$\varepsilon_{\rm CG}$ when using the new PCG iteration stopping rule.


\subsection{Evaluating the distance between convex polyhedra}
\label{NumTest2}
Let the two convex polyhedra ${\cal X}_1$ and ${\cal X}_2$ be described by 
the following two systems of linear inequalities: 
$$
{\cal X}_1 = \lbrace x_1: ~A_1^Tx_1\le b_1 \rbrace, \qquad
{\cal X}_2 = \lbrace x_2: ~A_2^Tx_2\le b_2 \rbrace,
$$ 
where $A_1\in R^{s\times n_1}$, ~$A_2\in R^{s\times n_2}$, and the vectors 
$x_1, x_2, b_1, b_2$ are of compatible dimensions. The original problem of evaluating the distance between ${\cal X}_1$ and ${\cal X}_2$ is (cf. \cite{Bo89}, where it was solved by the projected gradient method)
$$
x_*=\arg\min_{A_1^Tx_1\le b_1,A_2^Tx_2\le b_2}\|x_1-x_2\|^2/2,
\quad {\rm where} \quad x=[x_1^T,x_2^T]^T\in R^{2s}.
$$
We will use the following regularized/penalized approximate reformulation 
of the problem in terms of unconstrained convex piecewise quadratic minimization. Introducing the matrices 
$$
A = 
\left[
\begin{array}{cc}
A_1 &  0  \\
 0  & A_2
\end{array}
\right] \in R^{2s\times(n_1+n_2)},
\qquad 
B = 
\left[
\begin{array}{cc}
~I_s & -I_s  \\
-I_s &~~I_s
\end{array}
\right] \in R^{2s\times 2s},
$$
and the vector $b=[b_1^T,b_2^T]^T\in R^{n_1+n_2}$, we consider the problem
$$
x_*(\varepsilon) = \arg\min_{x\in R^{2s}} 
\left(\frac{\varepsilon}{2}\|x\|^2 + \frac{1}{2} x^TBx 
+ \frac{1}{2\varepsilon}\|(A^Tx-b)_+\|^2\right),
$$
where the regularization/penalty parameter $\varepsilon$ is a sufficiently small positive number (we have used $\varepsilon=10^{-4}$). The latter 
problem can readily be solved by adjusting the above described Algorithm~1  using $\delta=0$ and the following explicit expressions for the gradient 
and the generalized Hessian:
$$
g(x)=\varepsilon x + Bx +\varepsilon^{-1}A(A^Tx-b)_+, \qquad
H(x)=\varepsilon I + B + \varepsilon^{-1}A D(x) A^T,
$$
where $D(x)={\rm Diag}({\rm sign}(A^Tx-b)_+)$. When solving practical 
problems of evaluating the distance between two 3D convex polyhedra 
determined by their faces (so that $s=3$), the inexact Newton iterations 
are performed in $R^6$, and the cost of each iteration is proportional to 
the total number $n=n_1+n_2$ of the faces determining the two polyhedra. 
In this case, the explicit evaluation of $H(x)$ and the use of its Cholesky factorization is more preferable than the use of the CG method.

Test polyhedrons with $n/2$ faces each were centered at the points 
$e=[1,1,1]^T$ or $-e$ and defined as $A_1^T(x-e)\le b_1$ and 
$A_2^T(x+e)\le b_2$ with $b_1=b_2=[1,\ldots,1]\in R^{n/2}$, respectively.
The columns of matrices $A_1$ and $A_2$ were determined by $n/2$ 
quasirandom unit $3$-vectors generated with the use of logistic sequence 
(see, e.g.\cite{YBZS10} and references cited therein) 
$\xi_0=0.4$, ~$\xi_k=1-2\xi_{k-1}^2$. 
We used $A_1(i,j)=\xi_{20(i-1+3(j-1))}$ and similar for $A_2$; 
then the columns of these matrices were normalized to the unit length.
The corresponding performance results seem quite satisfactory, 
see Table~\ref{tab4}. 
Note that the lower bound $\|x_1-x_2\|\ge 2\sqrt{3}-2 \approx 1.464101$ 
always holds for the distance (since the two balls $\|x_1-e\|=1$ and 
$\|x_2+e\|=1$ are inscribed in the corresponding polyhedrons). 

\begin{table}
\begin{small}
\caption{Performance of the generalized Newton method 
($\varepsilon=10^{-4}$) for the problem of distance between 
two quasirandom convex polyhedrons with $n/2$ faces each} 
\label{tab4}
\begin{tabular}{p{1.8cm}p{1.9cm}p{1.9cm}p{1.9cm}p{1.9cm}p{1.6cm}}
\hline\noalign{\smallskip}
n & $\|x_1-x_2\|_2$ & $\|A^Tx_*-c\|_{\infty}$ & time(sec) & 
$\|g(x_*)\|_{\infty}$ & {\#}NewtIter\\ 
\noalign{\smallskip}\hline\noalign{\smallskip}
    8 & 0.001815 & 9.69--09 &$<0.001$ & 7.89--13 &  15 \\
   16 & 0.481528 & 8.63--05 &$<0.001$ & 1.27--13 &   3 \\
   32 & 0.795116 & 8.80--05 &$<0.001$ & 1.46--12 &  28 \\
   64 & 1.102286 & 1.32--04 &$<0.001$ & 5.58--13 &  13 \\
  128 & 1.446262 & 1.36--04 &$<0.001$ & 7.12--13 &  17 \\
  256 & 1.449913 & 9.54--05 &$<0.001$ & 4.37--13 &  11 \\
  512 & 1.460197 & 1.31--04 &$ 0.001$ & 8.16--13 &  15 \\
 1024 & 1.460063 & 1.46--04 &$ 0.002$ & 1.09--12 &  14 \\
 2048 & 1.463320 & 1.04--04 &$ 0.005$ & 6.58--13 &  19 \\
 4096 & 1.463766 & 1.26--04 &$ 0.009$ & 3.59--13 &  20 \\
 8192 & 1.463879 & 1.03--04 &$ 0.009$ & 8.32--14 &  12 \\
16384 & 1.463976 & 7.58--05 &$ 0.009$ & 1.64--12 &  13 \\
32768 & 1.464046 & 3.28--05 &$ 0.018$ & 1.54--12 &  13 \\
\noalign{\smallskip}\hline\noalign{\smallskip}
\end{tabular}
\end{small}
\end{table}

\begin{acknowledgement}
This work was supported by the Russian Foundation for Basic Research 
grant No.~17-07-00510 and by Program No.~26 of the Presidium of the 
Russian Academy of Sciences. The authors are grateful to the anonymous
referees and to Prof.~V.Garanzha for many useful comments which greatly
improved the exposition of the paper.

\end{acknowledgement}




\begin{thebibliography}{99.}

\bibitem{Ax76}
Axelsson, O.:
A class of iterative methods for finite element equations.
Computer Meth. Appl. Mech. Engrg. \textbf{9}, 123--137 (1976)


\bibitem {AK01}
Axelsson, O., Kaporin, I.E.:
Error norm estimation and stopping criteria in
preconditioned conjugate gradient iterations.
Numer. Linear Algebra Appls. \textbf{8} (4), 265--286 (2001)




\bibitem{Bo89}
Bobrow, J.E.: 
A direct minimization approach for obtaining the distance between convex polyhedra. The International Journal of Robotics Research, \textbf{8}(3), 65--76 (1989)


\bibitem{DS83}
Dembo R., Steihaug T.:
Truncated Newton algorithms for large-scale unconstrained optimization. 
Math. Program. \textbf{26}, 190--212 (1983)


\bibitem{DE03}
Dongarra, J., Eijkhout, V.:
Finite-choice algorithm optimization in Conjugate Gradients. 
Lapack Working Note 159, 
University of Tennessee Computer Science Report UT-CS-03-502 (2003)


\bibitem{GGE18}
Ganin, B.V., Golikov, A.I., Evtushenko, Y.G.:
Projective-dual method for solving systems of linear equations 
with nonnegative variables. 
Comput. Math. and Math. Phys. \textbf{58} (2), 159--169 (2018)


\bibitem{GK99}
Garanzha V.A., Kaporin I.E.:
Regularization of the barrier variational method of grid generation.
Comput. Math. and Math. Phys. \textbf{39} (9), 1426--1440 (1999)


\bibitem{GKK04} 
Garanzha, V., Kaporin, I., Konshin, I.:
Truncated Newton type solver with application to grid untangling problem.
Numer. Linear Algebra Appls. \textbf{11} (5-6), 525--533 (2004)


\bibitem{GGEN09}
Garanzha, V.A., Golikov, A.I., Evtushenko, Y.G., Nguen, M.K.:
Parallel implementation of Newton's method for solving 
large-scale linear programs. 
Comput. Math. and Math. Phys. \textbf{49} (8), 1303--1317 (2009)

\bibitem{HSN84}
Hiriart-Urruty, J. B., Strodiot, J. J., Nguyen, V. H.:  
Generalized Hessian matrix and second-order optimality conditions for problems with 
$C^{1,1}$ data. Applied mathematics and optimization, \textbf{11}(1), 43--56 (1984)

\bibitem{HS52}
Hestenes, M.R., Stiefel, E.L.: 
Methods of conjugate gradients for solving linear systems. 
J. Research Nat. Bur. Standards \textbf{49} (1), 409--436 (1952)


\bibitem {KA94}
Kaporin, I.E., Axelsson, O.:
On a class of nonlinear equation solvers based on the 
residual norm reduction over a sequence of affine subspaces.
SIAM J. Sci. Comput. \textbf{16} (1), 228--249 (1995)


\bibitem{Ka03} 
Kaporin, I.E.:
Using inner conjugate gradient iterations in solving 
large-scale sparse nonlinear optimization problems.
Comput. Math. and Math. Phys. \textbf{43} (6), 766--771 (2003)


\bibitem{KM11}
Kaporin, I.E., Milyukova, O.Y.:
The massively parallel preconditioned conjugate gradient method 
for the numerical solution of linear algebraic equations. 
In: Zhadan, V.G. (ed.) Collection of Papers of the Department 
of Applied Optimization of the Dorodnicyn Computing Center, 
pp. 132--157, Russian Academy of Sciences, Moscow (2011)


\bibitem{KMPN17}
Ketabchi, S., Moosaei, H., Parandegan, M.,  Navidi, H.:
Computing minimum norm solution of linear systems of equations by the 
generalized Newton method. 
Numerical Algebra, Control and Optimization, \textbf{7} (2), 113--119 (2017)


\bibitem{Ma02}
Mangasarian, O.L.: 
A finite Newton method for classification. 
Optimization Methods and Software, \textbf{17} (5), 913--929 (2002)


\bibitem{Ma04}
Mangasarian, O.L.: A Newton method for linear programming. 
Journal of Optimization Theory and Applications, \textbf{121} (1), 1--18 (2004)



\bibitem{YBZS10}
Yu, L., Barbot, J.P., Zheng, G., Sun, H.: 
Compressive sensing with chaotic sequence.
IEEE Signal Processing Letters, \textbf{17}(8), 731--734 (2010) 

\end{thebibliography}
\end{document}